\documentclass[a4paper,twoside,
11pt]{amsart}
\usepackage{amssymb,amsfonts,amsmath}
\begin{document}
\newtheorem{defin}{~~~~Definition}
\newtheorem{prop}{~~~~Proposition}[section]
\newtheorem{remark}{~~~~Remark}[section]
\newtheorem{cor}{~~~~Corollary}
\newtheorem{theor}{~~~~Theorem}
\newtheorem{lemma}{~~~~Lemma}[section]
\newtheorem{ass}{~~~~Assumption}
\newtheorem{con}{~~~~Conjecture}
\newtheorem{concl}{~~~~Conclusion}
\renewcommand{\theequation}{\thesection.\arabic{equation}}
\newcommand{\vf}{\varphi}
\title{Nurowski's conformal structures for
(2,5)-distributions via dynamics of abnormal extremals}
\author {Andrei Agrachev\address{ S.I.S.S.A., Via Beirut 2-4,
34014 Trieste Italy and Steklov Mathematical Institute,
ul.~Gubkina~8, 117966 Moscow Russia; email: agrachev@sissa.it} \and
Igor Zelenko \address{ S.I.S.S.A., Via Beirut 2-4, 34014 Trieste
Italy; email: zelenko@sissa.it}} 
\maketitle \markboth{Andrei Agrachev and Igor Zelenko} {Nurowski's
conformal structures for $(2,5)$-distributions}
\begin{abstract} As was shown recently by P. Nurowski,
to any  rank 2 maximally nonholonomic vector distribution on a
$5$-dimensional manifold $M$ one can assign the canonical conformal
structure of signature $(3,2)$.
His construction is based on the properties of the special
$12$-dimensional coframe bundle over $M$, which was distinguished by
E. Cartan during his famous construction of
the canonical coframe for this type of distributions on some
$14$-dimensional principal bundle over $M$.
The natural question is how "to see" the Nurowski conformal
structure of a $(2,5)$-distribution purely geometrically without the
preliminary construction of the canonical frame.
We
give rather simple answer to this question, using the notion of
abnormal extremals of $(2,5)$-distributions and the classical notion
of the osculating quadric for curves in the
projective plane. Our method is a particular case of a general
procedure for construction of algebra-geometric structures for a
wide class of distributions, which will be described elsewhere.
We also relate the fundamental
invariant of $(2,5)$-distribution, the Cartan covariant binary
biquadratic form, to the classical Wilczynski invariant of curves in
the projective plane.
\end{abstract}
\vskip .2in


\section{Introduction}
\indent

\subsection{Statement of the problem}
The following rather surprising fact was discovered by P. Nurowski
recently in \cite{nur}: any rank 2 maximally nonholonomic vector
distribution on a $5$-dimensional manifold $M$ possesses the
canonical conformal structure of signature $(3,2)$. His construction
is based on the Cartan famous paper \cite{cartan}, where  the
canonical coframe for this type of distributions on some
$14$-dimensional principal bundle $P$ over $M$ was constructed. In
the modern terminology Cartan constructed the Cartan normal
connection, valued on the split real form $\widetilde G_2$ of the
exceptional Lie algebra $G_2$, which is the algebra of infinitesimal
symmetries of the most symmetric distribution of the considered
type. It is well-known that the algebra $\tilde G_2$
can be naturally imbedded into the algebra $so(4,3)$, which in turn
is isomorphic to the conformal algebra of signature $(3,2)$.
Motivated by this fact, P. Nurowski noticed that a simple quadratic
expression, involving only the 1-forms from the mentioned canonical
coframe, which annihilate the fibers of $P$, is transformed
conformally along the fibers of $P$ and therefore defines the
canonical conformal structure on $M$. In his construction he uses
actually not the bundle $P$ but the special $12$-dimensional coframe
bundle, to which E. Cartan reduces the bundle of all coframes in the
first step of his reduction-prolongation procedure (see also section
4 below). Note that this reduction
is
highly nontrivial and purely algebraic in nature.
We asked ourselves, how "to see" the field of cones, defining the
Nurowski conformal structure of a $(2,5)$-distribution, purely
geometrically without exploiting the Cartan equivalence method, the
bundles of so high (with compare to $M$) dimensions
and
the properties of the group of symmetries of the most symmetric
model?
We give quite elementary answer to this question, using the dynamics
of so-called abnormal extremals (or singular curves) of
$(2,5)$-distributions, living in some $6$-dimensional bundle over
$M$ (so that the fibers are just one-dimensional), and the classical
notion of the osculating quadratic cones for curves in the
projective plane.


\subsection{The main constructions} We start our construction with
the description of characteristic curves of $(2,5)$-distribution
$D$. Let $D^l$ be the $l$-th power of the distribution $D$, i.e.,
$D^l=D^{l-1}+[D,D^{l-1}]$, $D^1=D$. The tuple $\bigl(\dim D(q),\dim
D^2(q),\dim D^3(q),\ldots\bigr)$ is called the \emph{small growth
vector of $D$ at $q$}. Throughout this paper we will consider
$(2,5)$-distributions with the small growth vector $(2,3,5)$ at any
point.
Let $\pi:T^*M\mapsto M$ be the canonical projection. For any
$\lambda\in T_q^*M$,
$q\in M$,
let $\varsigma_\lambda(\cdot)=\lambda(\pi_*\cdot)$ be the canonical
Liouville form and $\sigma=d\varsigma$ be the standard symplectic
structure on $T^*M$. Denote by $(D^l)^{\perp}\subset T^*M$ the
annihilator of the $l$th power $D^l$, namely
\begin{equation}
\label{annihil} (D^l)^{\perp}= \bigl\{\lambda\in T^*M:\,\,
\langle\lambda, v\rangle=0\,\,\forall\, v\in
D^l\bigl(\pi(\lambda)\bigr)\bigr\}.
\end{equation}
Since the submanifold $(D^2)^\perp$ has odd codimension in
$T^*M$, 
the kernels of the restriction
$\sigma|_{(D^2)^\perp}$ of $\sigma$ on $(D^2)^\perp$ are
not trivial. Moreover for
the points of
$(D^2)^\perp\backslash (D^3)^\perp$ these kernels are
one-dimensional.
They form the {\it
characteristic line distribution} of the distribution $D$ on
$(D^2)^\perp\backslash(D^3)^\perp$.
The integral curves of this line distribution define
 a {\it characteristic 1-foliation} of $D$
on $(D^2)^\perp\backslash(D^3)^\perp$.
 In Control Theory the leaves of this 1-foliation
are also called {\it regular abnormal extremals of $D$}.

Now let us define the same objects on the projectivization of
$(D^2)^\perp\backslash(D^3)^\perp$. First note that in the
considered case $(D^3)^\perp$ coincides with the zero section of
$T^*M$. So the projectivization of
$(D^2)^\perp\backslash(D^3)^\perp$ is just ${\mathbb P}(D^2)^\perp$.
Further note that the homotheties of the fibers of $(D^2)^\perp$
preserve the characteristic line distribution.
Therefore the projectivization induces on ${\mathbb P}(D^2)^\perp$
in a natural way the \emph { characteristic line distribution},
and the \emph
{characteristic $1$-foliation}, which will be denoted by
by ${\mathcal C}$ and $\bar{\mathcal C}$ respectively. The leaves of
the $1$-foliation $\bar {\mathcal C}$ will be called the \emph
{abnormal extremals of the distribution $D$} (on ${\mathbb
P}(D^2)^\perp$). The distribution ${\mathcal C}$ can be defined
equivalently in the following way: take  the corank $1$ distribution
on $(D^2)^\perp\backslash(D^3)^\perp$, given by
the Pfaffian equation 
$\varsigma|_{(D^2)^\perp}=0$ and push forward it under
projectivization to ${\mathbb P}(D^2)^\perp$. In this way we obtain
a corank 1 distribution on ${\mathbb P}(D^2)^\perp$, which will be
denoted by $\Delta$. The distribution $\Delta$ defines the
quasi-contact structure on the even dimensional manifold ${\mathbb
P}(D^2)^\perp$ and ${\mathcal C}$ is exactly the characteristic
distribution of this quasi-contact structure.

For the coordinate-free construction
it is more convenient to work with the
projectivization ${\mathbb P}(D^2)^\perp$ rather than with
$(D^2)^\perp$. For simplicity, the canonical projection from
${\mathbb P}(D^2)^\perp$ to $M$ will be denoted by  $\pi$, as in the
case of the cotangent bundle $T^*M$.
 Given a
segment $\gamma$ of abnormal extremal on  ${\mathbb
P}(D^2)^\perp$ one can construct two special (unparameterized) curves
in three and two-dimensional projective spaces respectively. For
this for any $\lambda\in {\mathbb P}(D^2)^\perp$ denote by
$V_\lambda$ the following subspace of $T_\lambda{\mathbb
P}(D^2)^\perp$
\begin{equation}
\label{prejac} V_\lambda= T_\lambda ({\mathbb
P}T^*_{\pi(\lambda)}M)\cap T_\lambda\mathbb P(D^2)^\perp.
 \end{equation}
In other words, $V_\lambda$ is  the tangent space  to the fiber
of ${\mathbb P}(D^2)^\perp$ at the point $\lambda$. Note that by
constructions $\dim\,V_\lambda=1$. Let $O_\gamma$ be a neighborhood
of $\gamma$ in ${\mathbb P}(D^2)^\perp$ such that
$N=O_\gamma /(\bar {\mathcal C }|_{O_\gamma})$
 is a
well-defined smooth manifold.
Let $\phi
:O_\gamma\to N$ be the canonical projection on the factor.

First note that $\widetilde \Delta \stackrel{def}{=}\phi_*\Delta$ is
well defined corank 1 distribution on $N$. Also the distribution
$\widetilde \Delta$ is contact distribution. Second let
\begin{equation}
\label{jacurve1}
{\mathcal E}_\gamma(\lambda)
{=} \phi_*\bigl(V_\lambda \bigr), \quad \forall
\lambda\in\gamma.
\end{equation}
The curve $\{{\mathcal E}_{\gamma}(\lambda), \lambda\in \gamma\}$ is
a curve in the three-dimensional projective space $\mathbb P
\widetilde \Delta (\gamma)$.

\begin{remark}
\label{jacrem} \rm { The curve $E_\gamma$ is closely related to the
so-called Jacobi curve $J_\gamma$ of the abnormal extremal
$\gamma$,
\begin{equation}
\label{prejac1} J_\gamma(\lambda)=\phi_*\Bigl(\{v\in
T_{\lambda}\mathbb P(D^2)^\perp:\,\pi_*\,v\in
D(\pi\bigl(\lambda)\bigr)\}\Bigr),
\end{equation}
which was  intensively used in our previous papers \cite{jac1},
\cite{zelvar},\cite{zelcart}, and \cite{zeldoub} for construction of
differential invariants and canonical frames of rank $2$
distributions. The curve $J_\gamma$ is the  curve of $2$-dimensional
subspaces in the $4$-dimensional space $\widetilde\Delta(\gamma)$.
Besides, since the distribution $\widetilde\Delta$ is contact, the
space $\widetilde D(\gamma)$ is endowed with the symplectic
structure, defined up to a constant factor. The curve $J_\gamma$ is
the curve of Lagrangian subspaces w.r.t. this symplectic form.
Finally, the relation between the curves $J_\gamma$ and ${\mathcal
E}_\gamma$ is given by the formula
\begin{equation}
 \label{JErel} J_\gamma(\lambda)=\pi_*\bigl([C,
V]_\lambda\}\bigr)={\rm span} \bigr\{\mathcal E_\gamma(\lambda),
\frac{d}{dt}\mathcal E_\gamma\bigr(\Gamma(t)\bigl)|_{t=0}\bigl\},
\end{equation}
 where $[C, V]_\lambda=\{[h, v](\lambda): h\in {\mathcal C}, v\in
V\,\,{\rm are}\,\, {\rm vector}\,\, {\rm fields}\}$ and
$\Gamma(\cdot)$ is a parameterization of $\gamma$ such that
$\Gamma(0)=0$. So, by (\ref{JErel}) the Jacobi curve can be uniquely
recovered from the curve ${\mathcal E}_\gamma$.} $\Box$
\end{remark}

Further fix some point $\lambda\in \gamma$ and denote by
$\Pi_\lambda:\widetilde\Delta(\gamma)\mapsto\widetilde
\Delta(\gamma)/{\mathcal E}_{\gamma}(\lambda)$ the canonical projection
to the factor-space.
Let $L_\lambda(\bar\lambda)
{=} \Pi_\lambda\bigl({\mathcal E}_\gamma(\bar\lambda)\bigr)$  for
any $\bar\lambda\in\gamma\backslash\{\lambda\}$. It is easy to show
(see formula (\ref{remov}) below) that one can define
$L_\lambda(\bar\lambda)$ also for $\bar\lambda=\lambda$ such that
for some neighborhood $I$ of $\lambda$ on $\gamma$ the curve $\{
L_\lambda(\bar\lambda), \bar\lambda\in I\}$ will be a smooth curve
in the projective plane ${\mathbb P}\Bigl(\
\bar\Delta(\gamma)/{\mathcal E}_\gamma(\lambda)\Bigr)$. Moreover,
it is easy to see (using, for instance,  Remark \ref{jacrem} above and
Propositions 3.5 from \cite{zelvar}) that this curve has no inflection points (or, equivalently, it is regular curve
in the terminology of section 2 below). In general,
let $Y$ be a three-dimensional linear space and $\xi$ be
a smooth curve without inflection points
in its projectivization. It is well known (see, for example,
\cite{cartanpr},
 p. 53 there ) that for any point $y\in \xi$ there exists the unique
quadric in ${\mathbb P}Y$, which has  contact of order four  with
$\xi$ at $y$, i.e. it has the same forth jet at $y$ as $\xi$. This
quadric is called the \emph {osculating quadric of $\xi$ at $y$}.
The two-dimensional quadratic cone in $Y$, corresponding to this
quadric, is called the \emph {osculating cone of $\xi$ at $y$}.

By construction, there is the  natural identification:
\begin{equation}
\label{natid0}
T_\gamma N
\sim
T_\lambda\mathbb P(D^2)^\perp/\mathcal C_\lambda
,
\end{equation}
which in turn implies that $T_\gamma N/{\mathcal E}_\gamma(\lambda)$
is naturally identified with $T_\lambda\mathbb P(D^2)^\perp/{\rm
span}\{\mathcal C_\lambda, V_\lambda\}$. On the other hand, since
$V_\lambda$ is the tangent space at $\lambda$ to the fiber of
${\mathbb P}(D^2)^\perp$, the latter space
is naturally identified with $T_{\pi(\lambda)}M/\pi_*{\mathcal
C}_\lambda$. So, finally one has the following natural
identification
\begin{equation}
\label{natid} T_\gamma N/{\mathcal E}_\gamma(\lambda)\sim
T_{\pi(\lambda)}M/\pi_*{\mathcal C}_\lambda.
\end{equation}
Then under the identification (\ref{natid}) one has
\begin{equation}
\label{natidd} \tilde \Delta(\gamma)/{\mathcal
E}_\gamma(\lambda)\sim \lambda^\perp/\pi_*{\mathcal C}_\lambda,
\end{equation}
where $\lambda^\perp=\{v\in T_qM: \langle\lambda, v\rangle=0\}.$

Now fix some point $q\in M$.  For any subspace $W$ of $T_qM$ denote
by ${\rm pr}_W:T_qM\mapsto T_qM/W$ the canonical projection and by
$\mathbb P(D^2)^\perp(q)$ the fiber of $\mathbb P(D^2)^\perp$ over
$q$. For any $\lambda\in \mathbb P(D^2)^\perp(q)$ take the
osculating cone  of the curve
$\{L_\lambda(\bar\lambda):\bar\lambda\in\gamma_\lambda\}$ at the
point $\lambda$, where $\gamma_\lambda$ is the abnormal extremal of
$D$, passing through $\lambda$. Under the identification
(\ref{natidd}) this osculating cone is the quadratic cone in the
three-dimensional linear space $\lambda^\perp/\pi_*{\mathcal
C}_\lambda$. Let ${\rm Con}(\lambda)$ be the preimage of this cone
under  ${\rm pr}_{\pi_*{\mathcal C}_\lambda}$. The set ${\rm
Con}(\lambda)$ is the degenerated 3-dimensional quadratic cone in
the four-dimensional linear subspace $\lambda^\perp$ of $T_q M$.
Finally, let $\Xi_q$ be the union of the cones ${\rm Con}(\lambda)$
for all $\lambda$ from the fiber $\mathbb P(D^2)^\perp(q)$,
\begin{equation}
\label{conus} \Xi_q=\bigcup_{\lambda\in \mathbb P(D^2)^\perp(q)}{\rm
Con}(\lambda).
\end{equation}
The set $\Xi_q$ is a four-dimensional cone in $T_qM$,
associated intrinsically with our
distribution. What is the structure
of this set?
First of all it is easy to see that $\Xi_q$ is an algebraic variety.
Indeed, in the appropriate coordinates the vector fields, which span
the line distributions ${\mathcal C}$ and $V_\lambda$, have the
components, which are polynomials on each fiber of ${\mathbb
P}(D^2)^\perp$. Therefore for any natural $k$ the mapping from a
fiber ${\mathbb P}(D^2)^\perp(q)$ to the space of $k$-jets of the
curves in the projective plane, which assigns to any $\lambda\in
{\mathbb P}(D^2)^\perp(q)$ the $k$-jet of the curve
$\{L_\lambda(\bar\lambda):\bar\lambda\in\gamma\}$,
is a rational  mapping. It implies without difficulties that the
quadratic cones ${\rm Con}(\lambda)$ are zero levels of quadratic
forms with coefficient, which are polynomials on the fiber ${\mathbb
P}(D^2)^\perp(q)$. Hence $\Xi_q$ is an algebraic variety.
The main results of the present paper are the following
two theorems

\begin{theor}
\label{theorem1}
 The set $\Xi_q$ is a quadratic cone of signature $(3,2)$ for any $q\in M$.
 \end{theor}
\begin{theor}
\label{theorem2} The conformal structure, defined by the field of
cones $\{\Xi_q\}_{q\in M}$, coincides with the Nurowski conformal
structure.
\end{theor}

Theorem \ref{theorem1} and \ref{theorem2} are  proved in sections
2 and 4 respectively.
In section 3, as a preparation to the proof of Theorem \ref{theorem2}, we
obtain the explicit formula for the constructed conformal structure
in terms of the structural functions of any frame naturally adapted
to the distribution. This formula is useful by itself.

{\subsection{The Cartan form of (2,5)-distributions and the
Wylczynski invariant of curves in the projective plane} Both the
discovery of the canonical structure for $(2,5)$-distribution in
\cite{nur} and the variational approach to the equivalence problem
of distributions via abnormal extremals and Jacobi curves, developed
in \cite{jac1} and \cite{zelvar}, gave the new geometric
interpretations  of the fundamental invariant of
$(2,5)$-distribution, the Cartan covariant binary biquadratic form
(an invariant homogeneous polynomial of degree 4 on each plane
$D(q)$). On one hand, it can be obtained from the Weyl tensor of the
canonical conformal structure. On the other hand, it can be
constructed from the special degree four differential on the Jacobi
curves, which are curves in Lagrange Grassmannians (\cite{zelcart}).
Using the constructions above, one can obtain  one more, totally
elementary interpretation
 of the  Cartan form of a $(2,5)$-distribution
in terms of the classical Wilczynski invariant (\cite{wil})
of the curves
$\{L_\lambda(\bar\lambda):\bar\lambda\in\gamma_\lambda\}$, which are
curves in the projective plane. Let us briefly describe how to do
it, referring to the Appendix for the facts from the theory of
curves in the projective spaces.
Actually by our method we can construct more than just
quadratic cone in the tangent space $T_qM$ of any point $q$. For this
take the union of the preimages  of the curves
$\{L_\lambda(\bar\lambda):\bar\lambda\in\gamma_\lambda\}$ under the
projection ${\rm pr}_{\pi_*{\mathcal C}_\lambda}$, where $\lambda$
runs over the fiber $\mathbb P (D^2)^\perp(q)$, namely, let
\begin{equation}
\label{buket} \Omega_q=\bigcup_{\lambda\in\mathbb P
(D^2)^\perp(q)}{\rm pr}_{\pi_*{\mathcal
C}_\lambda}^{-1}\Bigl(\{L_\lambda(\bar\lambda):\bar\lambda\in\gamma_\lambda\}
\Bigr).
\end{equation}
The set $\Omega_q$ is a conic hypersurface in $T_qM$.
In order to get the value of the Cartan form at some vector $v\in D(q)$ one can proceed as follows: First note that there exists exactly one point $\lambda$ on the fiber of ${\mathbb P}(D^2)^\perp$ over the point $q$ such that the projection to the base manifold $M$ of the abnormal extremal
$\gamma_\lambda$,  passing through $\lambda$, is tangent to $v$.
From (\ref{buket}) it follows that the intersection of the set $\Omega_q$ with the hyperplane
$\lambda^\perp$ is the hypersurface in $\lambda^\perp$, which is exactly the preimage of the curve $\{L_\lambda(\bar\lambda):\bar\lambda\in\gamma_\lambda\}$ under $pr_v$. Further recall that for any regular curve
in the projective plane one can construct the invariant degree $3$
differential, i.e. a degree $3$ homogeneous polynomial on the
tangent line to the curve at any point (the (first) Wilczinsky
invariant in the terminology of the Appendix below). Applying
Proposition \ref{4rem} from the Appendix to the curve
$\{L_\lambda(\bar\lambda):\bar\lambda\in\gamma_\lambda\}$ we get
that for this  curve
this degree 3 differential is equal to zero at the point $\lambda$.
Therefore the derivative of this differential at $\lambda$ is
well-defined degree $4$  homogeneous polynomial on the tangent space
to the curve
$\{L_\lambda(\bar\lambda):\bar\lambda\in\gamma_\lambda\}$ at
$\lambda$. Take a parameter $t$  on  $\gamma_\lambda$ such that
$\gamma_\lambda(0)=\lambda$ and $\pi_*\dot \gamma_\lambda(0)=v$. It
turns out the value of the Cartan form at $v$ is equal, up to the
constant factor $-\frac {1}{5}$,  to the value of the mention degree
$4$ homogeneous polynomial at $\frac{d}{dt}
L_\lambda(\gamma_\lambda(t))|_{t=0}$.
The last statement can be obtained immediately if one combines
Remark \ref{jacrem}, Proposition \ref{4rem} and Remark
\ref{willgrem} from the Appendix with Theorem 2 from \cite{zelcart}.
It also shows that the fundamental invariant of $(2,5)$-distribution
at a point $q$ can be obtained as an
invariant of the germ of the set $\Omega_q$ at $0$.

\subsection{Other types of distributions}

Now we say several words about similar constructions for other types
of distributions.
First note that shortly after Nurowski's discovery, R. Bryant gave a
construction of the canonical conformal structure of signature
$(3,3)$ for maximally nonholonomic $(3,6)$-distributions
(\cite{br36}), exploiting again the Cartan method of equivalence.
What are the analogs of the Nurowski and Bryant conformal structures
for other classes of distributions?
It turns out that the method of osculating cones in the projective
plane given in the present paper is a  particular case of rather
general procedure, which allows to assign natural algebra-geometric
structures to vector distributions from very wide class. Let us
describe this procedure very briefly. First, similarly to above one
can define the $1$-foliation of abnormal extremals of a distribution
$D$ in some submanifold ${\mathcal S}$ of the cotangent bundle,
where $\mathcal S=D^\perp$ in the case of odd rank and $\mathcal S$
is a codimension one subbundle of $D^\perp$ with algebraic fibers in
the case of even rank (in the case of rank 2 distributions $\mathcal
S=(D^2)^\perp$ )\footnote{ If rank of a distribution is even and not
$2$, then the abnormal extremals may not exist. However the
procedure can be extended also to this case after the
complexification of the fibers of $T^*M$.} . Second, studying the
dynamics of the fibers of $\mathcal S$ along the foliation of
abnormal extremals,
for any $q\in M$ one
gets an algebraic variety ${\mathfrak K}_q$ on the linear space
$T_qM\oplus D^\perp(q)$.
In the case of  maximally nonholonomic $(2,5)$ and $(3,6)$
distributions the projection of this variety  to $T_qM$ gives the
quadratic cone there, which defines the Nurowski and the Bryant
conformal structure respectively. Moreover, the variety ${\mathfrak
K}_q$ can be reconstructed from this projection. In general case,
the variety ${\mathfrak K}_q$ can not be recovered from its
projection  to $T_qM$ (for example, in many cases this projection is
the whole $T_qM$). However the constructed varieties   ${\mathfrak
K}_q$, $q\in M$ in the vector bundle $TM\oplus D^\perp$ contain the
fundamental information about the distribution.
The details of these constructions will be given in our forthcoming
paper \cite{alggeom}.

\section{ Proof of Theorem \ref{theorem1}}
\setcounter{equation}{0}
 \subsection{Preliminaries}
We start with several simple auxiliary facts
from geometry of
curves in projective spaces. Let $Z$ be an $k$-dimensional linear
space and $\zeta$ be a smooth curve in the projective space $\mathbb
P Z$ or, equivalently, a smooth curve of straight lines in $Z$. Take
some parameter $t\in \tilde I$ on $\zeta$,
where $\tilde I$ is some segment in $\mathbb R$. The curve
$\varepsilon:\tilde I\mapsto Z$ such that $\varepsilon(t)\in \xi(t)$
for all $t\in \tilde I$ will be called a \emph{representative of the
curve $\xi$, corresponding to the parameter $t$}.
A smooth curve $\zeta$ in the projective space $\mathbb P Z$ is
called \emph{regular}, if some (and therefore any)
its representative $\epsilon: \tilde I\mapsto Z$ satisfies
\begin{equation}
\label{spancond}
{\rm span} \{\varepsilon^{(i)}(t)\}_{i=0}^{k-1}
=Z,
\quad \forall t\in \tilde I,
\end{equation}
Assume that $\zeta$ is a regular curve and $t$ is some parameter on
it.
It is well-known that there exists the unique, up to the
multiplication on a constant, representative
$\varepsilon$,
corresponding to the parameter $t$
such that
\begin{equation}
\label{last0}
\varepsilon^{(k)}(t)=\sum_{i=0}^{k-2}B_i(t)\varepsilon^{(i)}(t),
\end{equation}
i.e. the coefficients near $\varepsilon^{(k-1)}(t)$ in the linear
decomposition of $\varepsilon^{(k)}(t)$ w.r.t. the basis
(\ref{spancond}) vanishes. Such representative of $\xi$ will be
called \emph{canonical w.r.t. the parameter $t$}.

\begin{remark}
\label{simprem}
{\rm
If in addition $\dim Z=k$ is even, $k=2m$, the space $Z$ is endowed with the
symplectic form $\bar\sigma$, and the $m$-dimensional subspace ${\rm span}
\{\varepsilon^{(i)}(t)\}_{i=0}^{m-1}$ is Lagrangian for any $t\in \tilde I$,
then the representative $t\mapsto \
\varepsilon(t)$ is canonical if and only if
\begin{equation}
\label{symp}
\bar\sigma
\bigl(\varepsilon^{(m)}(t),\varepsilon^{(m-1)}(t)\bigr)\equiv {\rm const}.
\end{equation}
Also in this case $B_{2m-3}(t)=B_{2m-2}'(t)$ and more general for
any $0\leq j\leq m-2$ the coefficient $B_{2j-1}$ is equal to some
polynomial expression w.r.t. the tuple of the coefficients
$\{B_{2l}\}_{l=j}^m$ and their derivatives. All statements in this remark can be easily
checked (see also \cite{zelprag}).} $\Box$
\end{remark}

Further, fix some point $\ell\in \zeta$, $\ell=\zeta(t_0)$. Note
that $\ell$ is a straight line in $Z$. Let $\Pi_\ell:Z\mapsto
Z/\ell$ be the canonical projection to the factor-space. If
$\varepsilon$ is a representative of  $\zeta$, corresponding to
the parameter $t$,
then the germ at $t_0$ of the following
map
\begin{equation}
\label{remov}
\tilde\varepsilon(t)=\left\{\begin{array}{ll}\Pi_\ell\bigl(\varepsilon(t)\bigr)& t\neq t_0\\
\Pi_\ell\bigl(\varepsilon'(t_0)\bigr)&t=t_0
\end{array}
\right.
\end{equation}
is the representative of the smooth curve in $(k-2)$-dimensional
projective space $\mathbb P(Z/\ell)$. This curve will be called the
\emph {reduction of the curve $\zeta$ by the point $\ell$}. The
point $y=\mathbb R \Pi_\ell\bigl(\varepsilon'(t_0)\bigr)$
will be called the \emph {point of the reduction, corresponding to
$\ell$}.

\begin{remark}
\label{regred}{\rm Note also that if the curve $\zeta$ is regular
then its reduction by $\ell$ is regular in a neighborhood of the
point
of the reduction, corresponding to $\ell$.} $\Box$
\end{remark}

Now let $Y$ be a three-dimensional linear space, $\xi$ be a regular curve in the projective plane $\mathbb P Y$,
 and $y$ be a point on $\xi$. In the following
lemma we give an explicit formula for the osculating quadric to
$\xi$ at $y$ in some coordinates under normalization
assumptions:
\begin{lemma}
\label{oscoord} Suppose that $(y_1,y_2)$ is a coordinate system in
$\mathbb P Y$ such that $y=(0,0)$ and $\tilde\varepsilon$ is a representative of $\xi$
such that $\tilde\varepsilon(0)=y$,
$\tilde\varepsilon(t)=\bigl(\alpha_1(t),\alpha_2(t)\bigr)$ in the chosen coordinates, and
\begin{equation}
\label{normass}
\dot\alpha_1(0)=\frac{1}{2},\,\,\dot\alpha_2(0)=\ddot\alpha_1(0)=0,
\ddot\alpha_2(0)=\frac{1}{3}.
\end{equation}
Then the osculating quadric to $\xi$ at $y$ has the following
equation
\begin{equation}
\label{osceq}
\frac{2}{3}y_1^2+2\alpha_2^{(3)}(0)y_1y_2-\left(4\alpha_1^{(3)}(0)+
6\bigl(\alpha_2^{(3)}(0)\bigr)^2-\frac{3}{2}\alpha_2^{(4)}(0)\right)y_2^2-y_2=0.
\end{equation}
\end{lemma}

The proof of Lemma \ref{oscoord} consists of straightforward
computations. The normalization assumptions (\ref{normass}) are
convenient when one makes a reduction.
Indeed, suppose that $\dim Z=4$, $Y=Z/\ell$ and
the curve
$\xi$ is the reduction of the curve $\zeta$ by a point $\ell$.
Let $\varepsilon$ be
a representative of $\zeta$ such that $\varepsilon(0)=\ell$, $\tilde\varepsilon$ be
the representative of $\xi$, satisfying (\ref{remov}) and $y=\tilde\varepsilon(0)$.
As a basis in $Y$ take the tuple $\Bigl(\Pi\bigl(\varepsilon'(0)\bigr),
\Pi\bigl(\varepsilon''(0)\bigr), \Pi\bigl(\varepsilon^{(3)}(0)\bigr)\Bigr)$. Let
$(Y_0, Y_1, Y_2)$ be the coordinates in $Y$ w.r.t. this basis. Then $(y_1, y_2)=
(Y_1/Y_0, Y_2/Y_0)$ are coordinates in a neighborhood of $y$ in $\mathbb P Y$.
If $\tilde\varepsilon(t)=\bigl(\alpha_1(t),\alpha_2(t)\bigr)$ in this coordinates,
then by direct check the functions $\alpha_1$ and $\alpha_2$
satisfy (\ref{normass}). Besides, if $\varepsilon$ is the canonical representative
of the curve $\zeta$, satisfying (\ref {last0}), where $k=4$, then
by simple computations
$$\alpha_1^{(3)}(0)=\frac{B_2(0)}{4},\quad \alpha_2^{(3)}(0)=0,\quad
\alpha_2^{(4)}(0)=\frac{B_2(0)}{5}.$$
Then from Lemma
\ref{oscoord} it follows that the osculating cone of $\xi$ at $y$ in $Y$
w.r.t. the chosen coordinates $(Y_0, Y_1, Y_2)$ has the following
equation
\begin{equation}
\label{osceq1}
\frac{2}{3}Y_1^2-\frac{7}{10}B_2(0) Y_2^2-Y_0Y_2=0.
\end{equation}
\begin{remark}
\label{rhorem} {\rm Comparing decomposition (\ref{last0}) and
formula (4.2) from \cite{zelvar}, we get that the function $B_2$
coincides, up to a constant factor, with the so-called Ricci
curvature $\rho$ of the parameterized curve $t\mapsto{\rm span}
\{\epsilon(t),\epsilon'(t)\}$ in the Grassmannian $G_2(Z)$ of
half-dimensional subspaces of the space $Z$. More precisely,
$B_2=-\frac{15}{2}\rho$. The definition of the Ricci curvature of
the curves in Grassmannians of half-dimensional subspaces of an
even-dimensional space can be found, for example, in subsection 2.3
of \cite{zelvar}.} $\Box$
\end{remark}
%
%

\subsection{Algebraic structure of the set $\Xi_q$} Note that the
curve $\mathcal E_{\gamma}$,
constructed  in the subsection 1.2  for a segment $\gamma$ of any
abnormal extremal, is a regular curve in the  corresponding
three-dimensional projective space (this fact follows, for example,
from Remark \ref{jacrem} above and Proposition 3.5 from
\cite{zelvar}). The curve $L_\lambda$, defined in the subsection 1.2
as well, is exactly the reduction of $\mathcal E_{\gamma}$ by
${\mathcal E}_\gamma(\lambda)$. Hence, by Remark \ref{regred} the
curve $L_\lambda$ is a regular curve in the corresponding projective
plane. Now we will use Lemma \ref{oscoord} for the explicit
calculation of the set $\Xi_q$ in appropriate coordinates. For this
it is more convenient to work with $(D^2)^\perp$ rather then with
its projectivization. To begin with given  two vector fields $X_1$,
$X_2$, constituting a local basis of distribution $D$,
one can construct a special vector field $\vec h_{_{X_1,X_2}}$
tangent to the characteristic $1$-foliation of $D$ on $(D^2)^\perp\backslash
(D^3)^\perp$. For this suppose
that
\begin{eqnarray}
&~&X_3=[X_1,X_2]\quad {\rm mod}\, D,\,\,\,
X_4=[X_1,[X_1,X_2]]=[X_1,X_3]\quad {\rm mod}\, D^2,
\label{x345} \\
&~&X_5=[X_2,[X_1,X_2]]=[X_2,X_3]\quad{\rm mod}\, D^2 \nonumber
\end{eqnarray}
The tuple $\{X_i\}_{i=1}^5$ is called the
\emph{adapted frame of the distribution $D$}.
Let us introduce ``quasi-impulses'' $u_i:T^*M\mapsto\mathbb R$,
$1\leq i\leq 5$, by the following formula:
$u_i(\tilde\lambda)=\langle\tilde\lambda,X_i(q)\rangle$.
For given function $G:T^*M\mapsto \mathbb R$ denote by $\vec G$ the
corresponding Hamiltonian vector field defined by the relation
$\sigma(\vec G,\cdot)=-d\,G(\cdot)$.
Then it is easy to show (see, for example \cite{zel}) that
the following vector field
\begin{equation}
\label{ham25} \vec h
_{_{X_1,X_2}}
=u_4 \vec u_2- u_5 \vec u_1
\end{equation}
is tangent to the characteristic $1$-foliation. In the sequel we
will work with the fixed local basis $(X_1, X_2)$, therefore we will
write $\vec h$ instead of $\vec h_{_{X_1,X_2}}$.

Further denote by ${\mathcal P}:T^*M\mapsto\mathbb P T^*M$ the
canonical projection. Let $\vec e$ be the Euler field, i.e., the
infinitesimal generator of homotheties on the fibers of $T^*M$.
Take a point $\tilde\lambda\in (D^2)^\perp\backslash (D^3)^\perp$
and let  $\tilde \gamma$ be a segment of the abnormal extremal of
$D$ in $(D^2)^\perp\backslash (D^3)^\perp$, passing through $\tilde
\lambda$. Let $\lambda={\mathcal P} (\tilde\lambda)$ and
$\gamma={\mathcal P}(\tilde\gamma)$.
Then, combining the identification (\ref{natid0}) and the
identification $T_
\lambda
\mathbb P(D^2)^\perp\sim
T_
{\tilde
\lambda
}
(D^2)^\perp/\mathbb R \vec e(\tilde \lambda)$, we get
\begin{equation}
\label{natid05}
T_\gamma N
\sim
T_{\tilde\lambda}\mathbb (D^2)^\perp/{\rm span}
\{\vec h(\tilde\lambda), \vec e(\tilde\lambda)\}.
\end{equation}
By analogy with (\ref{prejac}) let $\widetilde V_{\bar\lambda}$ be
the tangent to the fiber of $(D^2)^\perp$ at
$\bar\lambda\in(D^2)^\perp\backslash(D^3)^\perp$, namely,
\begin{equation}
\label{prejact} \widetilde V_{\bar\lambda}= T_{\bar\lambda}
(T^*_{\pi(\bar\lambda)}M)\cap T_{\bar\lambda}( D^2)^\perp.
 \end{equation}
Then  it is not hard to see that under identification
(\ref{natid05}) one has
\begin{equation}
\label{ident}
\mathcal E_{\gamma}\bigl(\mathcal P(e^{t H}\tilde \lambda)\bigr)
=(e^{-t
\vec h})_
*\bigl( \widetilde V (e^{t
\vec h}\tilde\lambda)\bigr)/
{\rm span}\{\vec h(\tilde \lambda),\vec e(\tilde\lambda)\}
\end{equation}
where $e^{t \vec h}$ is the flow generated by the vector field $
\vec h$.
The pair of functions $(u_4, u_5)$ defines the coordinates on
any fiber of $(D^2)^\perp$.
Consider the following vector field on $(D^2)^\perp\backslash(D^3)^\perp$, tangent
to the fibers of $(D^2)^\perp$:
\begin{equation}
\label{e125eq} \epsilon_1(\bar\lambda)=
\gamma_4(\bar\lambda) \partial_{u_4}+\gamma_5(\bar\lambda)
\partial_{u_5},
\end{equation}
where
\begin{equation}
\label{gammarel} \gamma_4 (\bar\lambda) u_5-\gamma_5(\bar\lambda)
u_4\equiv {\rm const}.
\end{equation}
Let\begin{equation}
\label{cane1J} v_{\tilde\lambda}(t)=(e^{-t\vec h})_*\epsilon_1\bigl(e^{t\vec
h}(\tilde\lambda)\bigr).
\end{equation}
Without introducing additional notations, we will look on
$v_{\tilde\lambda}(t)$ both as the element of $T_{\tilde\lambda} (D^2)^\perp$
and the element of the factor-space $T_{\tilde\lambda}\mathbb (D^2)^\perp/{\rm span}
\{\vec h(\tilde\lambda), \vec e(\tilde\lambda)\}$, i.e., the image of
$v_{\tilde\lambda}(t)$
under this factorization.
Combining Remarks \ref{jacrem}, \ref{simprem}, and the proof of Proposition 4.4 in \cite{zelvar},
we obtain  that under the identification (\ref{ident}) the curve
$t\mapsto v_{\tilde\lambda}(t)$ is the canonical representative of
the curve $\mathcal E_\gamma$ w.r.t. the parameterization
$t\mapsto\mathcal E_{\gamma}\bigl (\mathcal P(e^{t H}\tilde
\lambda)\bigr)$. Note that from (\ref{cane1J}) and the well known
property of Lie brackets it follows that
\begin{equation}
\label{lieb} v_{\tilde\lambda}^{(i)}(0)=({\rm ad} \,\vec
h)^i\epsilon_1(\tilde\lambda) \end{equation} (here, as usual, $({\rm
ad})\,\vec h {\mathcal X}=[\vec h, {\mathcal X}]$ for a vector field
$\mathcal X$). From (\ref{lieb}) and the definition of the canonical
representative (see (\ref{last0})) it follows that there exists the
tuple of functions $\{\mathcal B_i(\tilde\lambda)\}_{i=0}^2$ on
$(D^2)^\perp\backslash (D^3)^\perp$ such that
\begin{equation}
\label{ad4} ({\rm ad} \,\vec
h)^4\epsilon_1(\tilde\lambda)=\sum_{i=0}^2 {\mathcal
B}_i(\tilde\lambda)({\rm ad}\, \vec h)^i\epsilon_1(\tilde\lambda)
\end{equation}
(compare also with formula (4.10) in \cite{zelvar}). From the last
formula and (\ref{osceq1}) it follows that the set ${\rm
Con}(\mathcal P\tilde\lambda)$ defined in subsection 1.2
has the following description: If $(s,Y_0, Y_1, Y_2)$ are the
coordinates in $(\mathcal P\tilde\lambda)^\perp$ w.r.t. the basis
\begin{equation}
\label{basism} \Bigl(\pi_*\bigl(\vec h(\tilde\lambda)\bigr),
\pi_*\bigl([\vec h,\epsilon_1](\tilde\lambda)\bigr),\pi_*\Bigl(({\rm
ad}\, \vec h)^2\epsilon_1(\tilde\lambda)\Bigr),\pi_*\Bigl(({\rm
ad}\, \vec h)^3\epsilon_1(\tilde\lambda)\Bigr)\Bigr),
\end{equation}
 then the set ${\rm
Con}(\mathcal P\tilde\lambda)$ has the following equation in these
coordinates:
\begin{equation}
\label{conus1}
\frac{2}{3}Y_1^2-\frac{7}{10}{\mathcal B}_2(\tilde\lambda)
Y_2^2-Y_0Y_2=0.
\end{equation}
For simplicity take $\gamma_4=\frac{1}{u_5}$ and $\gamma_5\equiv 0$.
Then by direct calculations (one can use also formulas (4.24),
(4.35) and Lemmas 4.2 and 4.3 of \cite{zelcart}) it is not difficult
to show that
\begin{eqnarray}
&~&\pi_*\bigl([\vec h,\epsilon_1]\bigr)=-\frac{1}{u_5}X_2\quad  {\rm
mod}\{\mathbb R\pi_*\vec h\},\nonumber \\&~& \pi_*\bigl(({\rm ad}\,
\vec h)^2\epsilon_1\bigr)=\frac{l_1}{u_5} X_2+X_3 \quad  {\rm
mod}\{\mathbb R\pi_*\vec h\},\label {kisa}\\ &~&\pi_*\bigl(({\rm
ad}\, \vec h)^3\epsilon_1\bigr)=\frac{Q}{u_5}X_2+l_2X_3+u_4X_5-u_5
X_4 \quad {\rm mod}\{\mathbb R\pi_*\vec h\},\nonumber
\end{eqnarray}
where the functions $l_1$ and $l_2$ are linear
and the function $Q$ is a quadratic form
on each fiber of $(D^2)^\perp$ (in section 3 below we give the
explicit formulas for this functions in terms of the structural
functions of an adapted frame $\{X_i\}_{i=1}^5$, see formulas
(\ref{l1f})-(\ref{Qf})). Besides, from the proof of Theorem 3 in
\cite{zelvar} it follows that the function ${\mathcal B}_2$,
appearing in the equation (\ref{conus1}), is a quadratic form on
each fiber of $(D^2)^\perp(q)$ (see also explicit formula
(\ref{rhocar}) in section 3 below).

Now fixing the point $q\in M$, we will look at the functions $l_i$,
$Q$, and $\mathcal B_2$ as at the functions of $(u_4,u_5)$ (which
constitute coordinates on the fiber $(D^2)^\perp(q)$) and we will
write $l_i(u_4, u_5)$, $Q(u_4,u_5)$, and $\mathcal B_2(u_4,u_5)$
correspondingly. Then, running over the fiber $(D^2)^\perp(q)$ and
using (\ref{conus1}) and (\ref{kisa}), one can obtain by
straightforward calculations that if $(x_1,\ldots, x_5)$ are
coordinates in $T_qM$ w.r.t. the basis $\bigl(X_1(q),\ldots,
X_5(q)\bigr)$ then the set $\Xi_q$, defined in (\ref{conus}), has
the following equation in these coordinates:

\begin{equation}\begin{split}
\label{maineq} ~& x_1x_5-x_2x_4+\frac{2}{3}\Bigl(x_3-l_2(x_5,
-x_4)-\frac{3}{4}l_1(x_5,-x_4)\Bigr)^2-\\
~& \Bigl(\frac{7}{10}\mathcal B_2(x_5,-x_4)+Q(x_5,
-x_4)+\frac{3}{8}\bigl(l_1(x_5, -x_4)\bigr)^2\Bigr)=0.
\end{split}
\end{equation}
Taking into account the algebraic structure of the functions $l_1$,
$Q$, and ${\mathcal B}_2$,  we get that $\Xi_q$ is a quadratic cone.
Besides, it is easy to see that the quadratic form in the left
handside of (\ref{maineq}) has signature $(3,2)$. The proof of
Theorem 1 is completed.

The conformal structure, defined by the field of cones
$\{\Xi_q\}_{q\in M}$ will be called the \emph {canonical conformal
structure, associated with the distribution $D$}.

\section{The canonical conformal structure in terms of the structural
functions of an adapted frame}

\indent \setcounter{equation}{0}

Let, as before,  $\{X_i\}_{i=1}^5$, be an adapted frame of the
$(2,5)$-distribution $D$ and $u_i$, $1\leq i\leq 5$, be the
corresponding quasi-impulses. Denote by $c_{ji}^k$ the structural
functions of the frame $\{X_i\}_{i=1}^5$, i.e., the functions,
satisfying $$[X_i,X_j]=\sum_{k=1}^5 c_{ji}^k X_k.$$ We are going to
express the quadratic cone $\Xi_q$, $q\in M$, in terms of the
structural functions $c_{ji}^k$. For this we slightly modify the
computations produced in section 4 of \cite{zelcart}. As in
\cite{zelcart}, denote
\begin{eqnarray}
&~&\label{bmal}
 b(u_4, u_5)=
\frac{1}{3}\Bigl((c_{42}^4+c_{52}^5)u_4-(c_{41}^4+c_{51}^5)u_5
\Bigr),\\
&~&\label{b1def} b_1(u_4, u_5)=c_{32}^3u_4-c_{31}^3u_5,\\
&~&\label{aimal} \alpha_3(u_4,
u_5)=c_{52}^3u_4^2-(c_{42}^3+c_{51}^3)u_4u_5+
c_{41}^3 u_5^2,\\
&~&\label{Dbig}\Pi(u_4,
u_5)=(c_{32}^1+c_{53}^4)u_4^2+(c_{32}^2-c_{31}^1-c_{43}^4+c_{53}^5)u_4u_5-
(c_{31}^2+c_{43}^5) u_5^2.
\end{eqnarray}
Then from formulas (4.3) and  (4.36) of \cite{zelcart} it is not
hard to get that the function $l_1$, introduced in (\ref{kisa}),
satisfies
\begin{equation}
\label{l1f} l_1=b_1+3b
\end{equation}
(note that $\epsilon_1$ here is the same as $6\epsilon_1$ in
\cite{zelcart}). From formulas (4.3), (4.31)-(4.37), and (4.44) of
\cite{zelcart} one can obtain by direct computations that the
functions $l_2$ and $Q$ from (\ref{kisa}) satisfy:
\begin{eqnarray}
&~&\label{l2f} l_2=-3b,\\
&~&\label{Qf}Q=6\vec h(b)+\vec h(b_1)+\Pi-\alpha_3-b_1^2-3bb_1-9b^2.
\end{eqnarray}
Besides, combining Remark \ref{rhorem} with formula (4.56) from
\cite{zelcart}, we have
\begin{equation}
\label{rhocar} \mathcal B_2=2\alpha_3-\Pi-\vec h(b_1)-9\vec h(b)+
b_1^2+9b^2.
\end{equation}
 Finally, substituting expressions (\ref{l1f})-(\ref{rhocar}) into
(\ref{maineq}) we get easily that if $(x_1,\ldots, x_5)$ are
coordinates in $T_qM$ w.r.t. the basis $\bigl(X_1(q),\ldots,
X_5(q)\bigr)$ then the cone $\Xi_q$ has the following equation in
these coordinates
\begin{equation}\begin{split}
\label{maineq1} ~& x_1x_5-x_2x_4+\frac{2}{3}\Bigl(x_3-\frac{3}{4}(b_1-b)(x_5,-x_4)
\Bigr)^2-\\
~& \frac{3}{10}\Bigl(\Pi+\frac{4}{3}\alpha_3+\vec h(b_1-b)
+\frac{1}{4}(b_1-b)(b_1+9b)\Bigr)(x_5, -x_4)=0.
\end{split}
\end{equation}
In the last formula for shortness, instead of writing the arguments
$(x_5, -x_4)$ after each function $b$, $b_1$, $\alpha_3$ and $\Pi$,
we write them after the polynomial expressions, involving these
functions. So, for example, $(b_1-b)(x_5,-x_4)$ means
$b_1(x_5,-x_4)-b(x_5,-x_4)$. Note that from
(\ref{bmal})-(\ref{Dbig}) and (\ref{ham25}) it follows that
 $b$ and $b_1$ are linear forms in $u_4$ and $u_5$, while $\alpha_3$, $\Pi$, $\vec h(b_1-b)$
 are quadratic forms in $u_4$ and $u_5$. Therefore the equation
 (\ref{maineq1}) together with expressions (\ref{bmal})-(\ref{Dbig}) and
 (\ref{ham25}) gives the
explicit formula for the cone $\Xi_q$ it terms of the structural
functions of the adapted frame $\{X_i\}_{i=1}^5$.

\section{
Proof of Theorem 2} \indent \setcounter{equation}{0}

We start with the brief description of the construction of the
Nurowski conformal structure, associated with the
$(2,5)$-distribution $D$. His construction is based strongly on the
chapter VI of the E. Cartan paper \cite{cartan}, where the canonical
coframe for this type of distributions on some $14$-dimensional
principal bundle over $M$ was constructed.
%
To begin with let $\{\omega_i\}_{i=1}^5$
be a coframe on $M$
such that the rank $2$ distribution $D$ is the annihilator
of the first three elements of this coframe, namely,

\begin{equation}
\label{annih}
 D(q)=\{v\in T_qM;\,\,
\omega_1(v)=\omega_2(v)=\omega_3(v)=0\},\quad \forall q\in M
\end{equation}
Obviously the bundle of all coframes satisfying (\ref{annih}) has
19-dimensional fibers.
After quite cumbersome algebraic manipulation, E. Cartan succeed to
distinguish a special
$12$-dimensional subbundle ${\mathcal K}(M)$ of this coframe bundle
(having $7$-dimensional fibers), which is characterized as follows:
A coframe $\{\omega_i\}_{i=1}^5$ belongs to the coframe bundle
${\mathcal K}(M)$
if and only if it satisfies the following structural equations
(formula (VI.5) in \cite{cartan}):

\begin{eqnarray}
&~&d\,\omega_1=\omega_1\wedge(2\overline\omega_1+\overline\omega_4)+
\omega_2\wedge \overline\omega_2+\omega_3\wedge\omega_4\nonumber\\
&~&d\,\omega_2=\omega_1\wedge\overline\omega_3+\omega_2\wedge (
\overline\omega_1+2\overline\omega_4)+\omega_3\wedge\omega_5\nonumber\\
&~&d\,\omega_3=\omega_1\wedge\overline\omega_5+
\omega_2\wedge\overline\omega_6+\omega_3\wedge(\overline\omega_1
+\overline\omega_4)+\omega_4\wedge\omega_5\label{costruct}\\
&~&d\,\omega_4=\omega_1\wedge\overline\omega_7+
\frac{4}{3}\omega_3\wedge\overline\omega_6+
\omega_4\wedge\overline\omega_1+\omega_5\wedge\overline\omega_2\nonumber\\
&~&d\,\omega_5=\omega_2\wedge\overline\omega_7-
\frac{4}{3}\omega_3\wedge\overline\omega_5+
\omega_4\wedge\overline\omega_3+\omega_5\wedge\overline\omega_4,\nonumber
\end{eqnarray}
where $\overline \omega_j$, $1\leq j\leq 7$, are new $1$-forms. It
turns out that a further "reduction" of the coframe bundle
${\mathcal K}_1(M)$ (i.e., a selection of an intrinsic proper
subbundle of it) is impossible. Of course, it is only the first step
in Cartan's construction of the canonical coframe, but this step is
already enough for our purposes.
The remarkable observation of P. Nurowski was  that the field of cones
$\{{\mathcal N}_q\}_{q\in M}$, satisfying
\begin{equation}
\label{nurcone} \mathcal N_q=\left\{X\in
T_qM:\omega_1(X)\omega_5(X)-\omega_2(X)\omega_4(X)+
\frac{2}{3}\bigl(\omega_3(X)\bigr)^2=0
\right\},
\end{equation}
does not depend on the choice of a coframe $\{\omega_i\}_{i=1}^5$,
satisfying (\ref{costruct}).

Now we have all tools to prove that our cones $\Xi_q$ coincides with
the Nurowski cones ${\mathcal N}_q$. Fix some coframe
$\{\omega_i\}_{i=1}^5$ satisfying (\ref{costruct}).
Let $\{\widetilde X_k\}_{k=1}^5$, be the frame on $M$ dual to the
coframe $\{\omega_i\}_{i=1}^5$, namely, 
$\omega_i(\widetilde X_k)=\delta_{i,k}$.
Denote by \begin{equation} \label{XY} X_k=\widetilde X_{6-k},\quad
1\leq k\leq 5 .\end{equation} Any frame $\{X_i\}_{i=1}^5$ obtained
in such way will be called the {\it Cartan frame} of the
distribution $D$. From (\ref{costruct}) it is not difficult to show
that any Cartan frame is adapted frame of the distribution $D$ and
to find its structural functions (see formulas (5.5)-(5.13) in
\cite{zelcart}). Moreover, as already was noticed in \cite{zelcart}
(formulas (5.15) and (5.16) there), for any Cartan's frame
\begin{equation}
\label{cartid} b_1=b,\quad \Pi=-\frac{4}{3}\alpha_3.
\end{equation}
Substituting the last relations into (\ref{maineq1}) we obtain that
for any Cartan frame $\{X_i\}_{i=1}^5$, if $(x_1,\ldots, x_5)$ are
coordinates in $T_qM$ w.r.t. the basis $\bigl(X_1(q),\ldots,
X_5(q)\bigr)$ then the cone $\Xi_q$ has the following equation in
these coordinates
\begin{equation}
\label{maineq2} x_1x_5-x_2x_4+\frac{2}{3}x_3^2=0.
\end{equation}
Comparing this equation with (\ref{nurcone})  and taking into
account that the coframe $\{\omega_i\}_{i=1}^5$ is dual to the frame
$\{X_{6-i}\}_{i=1}$, we get immediately that $\Xi_q={\mathcal N}_q$,
which completes the proof of Theorem 2.

Actually one can relate to Nurowski's equation (\ref{nurcone}) for
the canonical cones, as to the particular case of our equation
(\ref{maineq1}) applied to a Cartan frame.

\section{Appendix}

In  subsection 1.3 we gave the description of the Cartan covariant
binary biquadratic form of $(2,5)$-distribution in terms of the
Wilczynski invariant of certain curves in the projective plane. In
the present Appendix we collect several facts about the invariants
of curves in the projective spaces and their reductions, on which
this description was based.

First of all let us recall, what the Wilczynski invariants of curves
in projective spaces (\cite{wil},\cite{doub}) are.
As in subsection 2.1 , let $Z$ be an $k$-dimensional linear space
and $\zeta$ be a regular curve in the projective space $\mathbb P
Z$. It turns out that the set of all parameters $t$ on $\zeta$, such
that for their canonical representative $\varepsilon$ the
coefficient $B_{k-2}$ in the decomposition (\ref{last0}) is
identically equal to zero, defines {\emph the canonical projective
structure on the curve $\zeta$} (i.e., any two parameters from this
set are transformed one to another by M\"{o}bius transformation). It
follows from the following reparameterization rule for the
coefficient $B_{k-2}$: if $\tau$ is another parameter,
$t=\varphi(\tau)$, and $\widetilde B_{k-2}$ is the coefficient from
the decomposition (\ref{last0}) corresponding to the canonical
representative of $\zeta$ w.r.t. the parameter $\tau$, then
\begin{equation}
 \label{rhorep}
 \widetilde
 B_{k-2}(\tau)=\dot\vf(\tau)^2B_{k-2}(\vf(\tau))-\frac{n(n^2-1)}{6}
\mathbb{S}(\vf),
\end{equation}
 where $\mathbb{S}(\vf)$ is a
Schwarzian derivative of $\vf$,
$\mathbb S(\varphi)=
=\frac {d}{dt}\Bigl(\frac {\varphi''}{2\,\varphi'}\Bigr)
-\Bigl(\frac{\varphi''}{2\,\varphi'}\Bigr)^2$. Any parameter from
the canonical projective structure of $\zeta$ is called \emph{ a
projective parameter}.

Now let $t$ be a projective parameter on $\zeta$. E. Wilczynski
showed that for any $i$, $1\leq i\leq k-2$, the following degree
$i+2$ differentials
\begin{equation}
\label{willy} {\mathcal
W}_i(t)\stackrel{def}{=}\frac{(i+1)!}{(2i+2)!}\left(\sum_{j=1}^{i}(-1)^{j-1}
\frac{
(2i-j+3)!(k-i+j-3)!}{(i+2-j)!j!}B_{k-3-i+j}^{(j-1)}(t)\right)(dt)^{i+2}
\end{equation}
 on
$\zeta$ does not depend on the choice of the projective parameters a
The form ${\mathcal W}_i$ is called \emph{the $i$-th Wilczynski
invariant of the curve $\zeta$}. In particular, the first and the
second Wilczynski invariants, which will be used in the sequel, have
the form

\begin{eqnarray}
&~&\label{willy1} \mathcal W_1=B_{k-3}(t)(dt)^3,\quad k\geq 3; \\
&~&\label{willy2} \mathcal
W_2=\left(B_{k-4}(t)-\frac{k-3}{4}B_{k-3}'(t)\right)(dt)^4,\quad k\geq 4
\end{eqnarray}

\begin{remark}
\label{willgrem} {\rm In a continuation of Remark \ref{simprem}, if
 $\dim Z=k$ is even, $k=2m$, the space $Z$ is endowed with the
symplectic form $\bar\sigma$,  the $m$-dimensional subspaces ${\rm span}
\{\varepsilon^{(i)}(t)\}_{i=0}^{m-1}$ is Lagrangian for any $t$, and
$t$ is a projective parameter, then
$B_{2m-3}\equiv B_{2m-2}'\equiv 0$. So,
\begin{equation}
\label{willg} {\mathcal W}_1\equiv 0, \quad \mathcal
W_2=B_{2m-4}(t)dt^4.
\end{equation}

Besides, it can be shown that ${\mathcal W}_2$ coincides, up to some
constant factor $C_m$, with the fundamental form of the curve of
Lagrangian subspaces ${\rm span}
\{\varepsilon^{(i)}(t)\}_{i=0}^{m-1}$, introduced in \cite{jac1}
(for its definition see also \cite[subsection 2.3]{zelvar}). In the
case $m=2$, corresponding to $(2,5)$-distributions, this fact
follows from formula (4.2) of \cite{zelvar} and the factor $C_2$ is
equal to $35$. }
$\Box$
\end{remark}

Now let $Y$ be a three-dimensional linear space, $\xi$ is a regular
curve in the projective space $\mathbb P Y$ and the first Wilczynski
invariant $\mathcal W_1$ of $\xi$ is equal to zero at some point
$y\in\xi$. If $t$ is a projective parameter, $t=t_0$ corresponds to
the point $y$ and $B_{k-3}(t)$ is the corresponding coefficient in
the decomposition (\ref{last0}), then $B_{k-3}(t_0)=0$ and the
following form $\mathcal W_1^{(1)}=B_{k-3}'(t_0)dt^4$ is the
well-defined degree $4$ differential at the point $y$ on the curve.
This degree $4$ differential will be called the \emph{derivative of
the first Wilczynski invariant at $\ell$}. Note also that in the
considered case there is only one  Wilczynski invariant, the first
one, so, referring to it, one usually  omits the word "first". Using
formulas (\ref{willg}),
it is not difficult to get the following

\begin{prop}
\label{4rem} Suppose that
 $\dim Z=4$, the space $Z$ is endowed with the
symplectic form $\bar\sigma$,  $\zeta$ is a regular curve in
$\mathbb P Z$ such that for its representative $t\mapsto
\varepsilon(t)$ all $2$-dimensional subspaces ${\rm span}
\{\varepsilon(t),\varepsilon'(t)\}$ are Lagrangian. If $\xi$ is the
reduction of the curve $\zeta$ by a  point $\ell$, then the
Wilczynski invariant of $\xi$ at the point $y$ of the reduction,
corresponding to  $\ell$,
is equal to zero and the derivative of the  Wilczynski invariant of
the curve $\xi$ at $y$ is equal, up to the factor $\frac{1}{5}$, to
the second Wilczynski invariant of the curve $\zeta$ at $\ell$.
\end{prop}

{\bf Acknowledgements} We thank Pawel Nurowski for very inspiring
course of lectures and discussions during his visit of SISSA in
January-February 2005, which stimulate us to work on this problem.

\end{document}